\documentclass{article}
\usepackage{subfig, enumerate}
\usepackage[section] {placeins}
\usepackage{amsthm}
\usepackage{amssymb,amsmath, mathdots}
\usepackage{anysize}
\usepackage{psfrag}
\usepackage{url}
\usepackage[normalem]{ulem}
\newtheorem{theorem}{Theorem}[section]
\newtheorem{definition}{Definition}[section]
\newtheorem{ex}{Example}[section]

\newtheorem{conj}{Conjecture}[section]

\numberwithin{table}{section}
\begin{document}

\title{On the orbits associated with the Collatz conjecture.}

\author{Louis H. Kauffman\\
Department of Mathematics, Statistics and Computer Science \\ 851 South Morgan Street   \\ University of Illinois at Chicago\\
Chicago, Illinois 60607-7045 USA\\ and\\ Department of Mechanics and Mathematics\\ Novosibirsk State University\\Novosibirsk, Russia \\
\texttt{kauffman@uic.edu}\\
and\\
Pedro Lopes\\
Center for Mathematical Analysis, Geometry, and Dynamical Systems, \\
Department of Mathematics, \\
Instituto Superior T\'{e}cnico, Universidade de Lisboa\\
1049-001 Lisbon, Portugal \\
\texttt{pelopes@math.tecnico.ulisboa.pt}}
%\date{September 18, 2018}
\maketitle
\qquad \qquad \qquad \qquad \qquad \qquad To Jos\'e Sousa Ramos, {\it in memoriam}
\bigbreak
\begin{abstract}
This article is based  upon previous work by Sousa Ramos and his collaborators. They first prove that the existence of only one orbit associated with the Collatz conjecture is equivalent to the determinant of each matrix of a certain sequence of matrices to have the same value. These matrices are called Collatz matrices. The second step in their work would be to calculate this determinant for each of the Collatz matrices. Having calculated this determinant for the first few terms of the sequence of matrices, their plan was to prove the determinant of the current term equals the determinant of the previous one. Unfortunately,  they could not prove it for the cases where the dimensions of the matrices are 26+54l or 44+54l, where l is a positive integer. In the current article we improve on these results.
\end{abstract}

Keywords: Collatz conjecture; Recurrence; Determinants; Permutations

Mathematics Subject Classification 2010: 14G10; 11D04; 15A15; 11D79; 11B83

\section{Introduction}\label{sec:intro}

The Collatz conjecture is a well-known conjecture about the asymptotic behaviour by the iterates of a certain function. It expects that these iterates eventually lie on a unique orbit, no matter which initial input is chosen. In the current article we improve the results concerning the  uniqueness of this orbit  obtained  in \cite{JFA}.

This is the version of the Collatz function we work with: $$f(n)=\begin{cases}
\frac{3n+1}{2},& \quad \text{ if $n$ is odd,}\\
\frac{n}{2},& \quad \text{ if $n$ is even.}
\end{cases}$$ Let $${\cal O}_1=\{ 1, 2 \} = {\cal O}_2.$$ We note that for $n=1$ (respect., $n=2$) the iterates all lie in ${\cal O}_1$. For $n=3$, the sequence of iterates begins with $$5, 8, 4, 2, 1, \dots $$ and the iterates eventually all lie in ${\cal O}_1$. We would like to prove that ${\cal O}_1$ is the only eventual orbit - or periodic orbit in the terminology of \cite{JFA}. We next present the notation, terminology and results in \cite{JFA} in order to position our own results.

\begin{definition}[Collatz Matrix]\label{def:CollatzM}
For each integer $k$ greater than $1$, we define the Collatz Matrix, denoted $M_k$, in the following way:
\begin{enumerate}[(i)]
\item It is a square $k\times k$ matrix;
\item Each of its diagonal entries is $1$;
\item For each even $1<i\leq \lfloor k/2\rfloor$, the $(i,i/2)$-entry is $x$;
\item For each odd $i$ such that $\frac{3i+1}{2}\leq k$, the $(i,\frac{3i+1}{2})$-entry is $x$;
\item Any other entry is $0$.
\end{enumerate}
\end{definition}

\begin{ex}\label{ex:CollatzM}
The following are examples of Collatz Matrices:
\begin{enumerate}[(i)]
\item
\begin{equation*}
M_2=
\begin{pmatrix}
1 & x\\
x & 1
\end{pmatrix}
\end{equation*}

\item
\begin{equation*}
M_3=
\begin{pmatrix}
1 & x & 0\\
x & 1 & 0\\
0 & 0 & 1
\end{pmatrix}
\end{equation*}

\item
\begin{equation*}
M_4=
\begin{pmatrix}
1 & x & 0 & 0\\
x & 1 & 0 & 0\\
0 & 0 & 1 & 0\\
0 & x & 0 & 1
\end{pmatrix}
\end{equation*}

\item
\begin{equation*}
M_5=
\begin{pmatrix}
1 & x & 0 & 0 & 0\\
x & 1 & 0 & 0 & 0\\
0 & 0 & 1 & 0 & x\\
0 & x & 0 & 1 & 0\\
0 & 0 & 0 & 0 & 1
\end{pmatrix}
\end{equation*}
\end{enumerate}
\end{ex}

\begin{conj}[Collatz Orbit Conjecture]\label{conj:Collatz}
${\cal O}_1$ is the only eventual finite orbit that the Collatz function admits.
\end{conj}

\begin{theorem}[\cite{JFA}]\label{thm:JFA}
For any integer $k>1$, $\det M_k = 1-x^2$ is equivalent to  the veracity of the Collatz Orbit Conjecture.
\end{theorem}

We note that, for $2\leq k \leq 5$, $$\det M_k = 1-x^2 .$$ The strategy in \cite{JFA} was, then, to prove that for any positive integer $k$ greater than $2$, $$\det M_{k}=\det M_{k-1} .$$ Unfortunately, they were not able to prove it for $$k=26+54l \qquad \text{ or } \qquad k=44+54l ,$$   where $l$ is a positive integer.

It should be noted that, for $k\neq 8+18l$, the computation is fairly straightforward, based on inspection of the last column  or the last row of $M_k$. For odd $k$, the last row only has one  non-zero entry, a $1$ in the $(k,k)$ position. Thus, by Laplace expansion over the last row, the recurrence relation is obtained for  odd $k$. Now, for even $k$. For some of these cases, there are two elements in the last row but there is only one in the last column, the $1$ in the diagonal entry. Here we use Laplace expansion along the last column. There remain the cases where there are two elements both in the last row and in the last column.  The elements in the last column are the $1$ in the $(k,k)$ position and the  $x$ in the $(\frac{2k-1}{3}, k)$ position. Laplace expansion over the last column is applied again but now we would like to prove that the minor matrix associated with the $x$ entry has zero determinant. Keeping to the approach used in \cite{JFA}, instead of this minor matrix we consider the matrix $$\widetilde{M}_{k-1}$$ which is a Collatz matrix except for the $\frac{2k-1}{3}$ row: it has an $x$ for the $(\frac{2k-1}{3}, \frac{k}{2})$ entry and all other entries (in this row) are zero. Specifically, this $\widetilde{M}_{k-1}$ is obtained from the minor matrix referred to above by a cyclic permutation of its last rows: row $k-1$ becomes row $\frac{2k-1}{3}$, which becomes row $\frac{2k-1}{3} +1$, and so on. Resuming the narrative, proving $\det \widetilde{M}_{k-1}=0$ is fairly  straightforward except for $k=18+18l$. For $k=8+54l$, Sousa Ramos and collaborators were able to prove the existence of a non-trivial solution of the system of linear homogeneous equations whose matrix of coefficients is $\widetilde{M}_{k-1}$. As remarked above, it remains to prove $\det \widetilde{M}_{k-1}=0$ for $k=26+54l$ and for $k=44+54l$. In the current article we are able to clear infinitely many of these cases - but not all. Also, we believe that our methodology provides a simpler solution for the case $k=8+54l$.
\subsection{Organization}\label{subsec:org}
The results are presented in Section \ref{sec:results} and the proofs are presented in Section \ref{sec:calcs}.
\subsection{Acknowledgements}\label{subsec:ack}

Kauffman's work was supported by the Laboratory of Topology and Dynamics, Novosibirsk State University (contract no. 14.Y26.31.0025 with the Ministry of Education and Science of the Russian Federation).

Lopes acknowledges support from FCT (Funda\c c\~ao para a Ci\^encia e a Tecnologia), Portugal, through project FCT PTDC/MAT-PUR/31089/2017, ``Higher Structures and Applications''.

\section{Results}\label{sec:results}

Theorem \ref{thm:results} below states our results; they supplement the results obtained in \cite{JFA}.

\begin{theorem}\label{thm:results}
Let $k$ be a positive integer. Let $M_k$ be a Collatz matrix.
\begin{enumerate}[(a)]
\item Assume further $k=44+54l$ for some positive integer $l$. Then, $\det M_k = \det M_{k-1}$, in the following instances:
\begin{enumerate}[(i)]
\item If $3\,|\,l$ or $3\,|\,(l-1)$.
\item Or, if $l=2+3l_1$, and $3\,|\, (l_1-1)$.
\item Or, if $l=2+3l_1$, and
\begin{enumerate}[(1)]
\item $l_1=3l_2$, and $3\,|\, l_2$; \qquad or
\item if $l_1=2+3l_2$, and $\bigg(3\,|\, l_2 \text{ or } 3\,|\,(l_2-1)\bigg)$.
\end{enumerate}
\end{enumerate}

\item Now assume $k=26+54l$ for some positive integer $l$. Then, $\det M_k = \det M_{k-1}$, in the following instances:
\begin{enumerate}[(i)]
\item If $3\,|\,(l-2)$.
\item Or,
\begin{enumerate}[(1)]
\item if $l=1+3l_1$, and $3\,|\, l_1$.
\item if $l=3l_1$, and $3\,|\, l_1$.
\end{enumerate}
\end{enumerate}
\end{enumerate}
\end{theorem}

{\bf Remark} As the reader can see in detail in Section \ref{sec:calcs}, we find that to show that  $\det M_k = \det M_{k-1}$ it is sufficient to show that k/2 does not lie on a specific closed (inverse) Collatz orbit.  This is the basis for the work in proving Theorem \ref{thm:results}.

\section{Calculations and Proofs}\label{sec:calcs}

We would like to prove that the sequence of Collatz matrices $(M_j)$ satisfies the property $$\det M_k = \det M_{k-1}.$$ Knowing that for the small values of $k$, $\det M_k=1-x^2$, this would imply that $\det M_k=1-x^2$, for all $k$, which would further imply that there is only one orbit associated with the Collatz function, that which contains the number $1$ (\cite{JFA}). We note that in \cite{JFA} much work has already been done in this direction, there remaining to be proved that $\det M_k = \det M_{k-1}$ only for $k$'s of the sort: $$k=44+54l \qquad \text{ or } \qquad k=26+54l \qquad \qquad \text{ where $l$ is a positive integer}.$$ For this sort of $k$, the last column of the Collatz matrix possesses two non-null entries.  A $1$ in the last row and an $x$ in row $(2k-1)/3$. Upon Laplace expansion about this last column we obtain $$\det M_k = 1\cdot \det M_{k-1}+x(-1)^{k+\frac{2k-1}{3}}\det M'_{k-1} .$$ We would then like to prove that $\det M'_{k-1}=0$. Instead of $M'_{k-1}$ we will work with a matrix obtained from this one by cyclic permutation of its last rows, namely row $k-1$ goes over to row $(2k-1)/3$ which goes over to row $(2k-1)/3+1$, and so on. We denote this matrix $$\widetilde{M}_{k-1} .$$ Here is a concrete example for $k=8$: \begin{equation*}
M_{8}=
\begin{pmatrix}
1 & x & 0 & 0 & 0 & 0 & 0 & 0\\
x & 1 & 0 & 0 & 0 & 0 & 0 & 0\\
0 & 0 & 1 & 0 & x & 0 & 0 & 0\\
0 & x & 0 & 1 & 0 & 0 & 0 & 0\\
0 & 0 & 0 & 0 & 1 & 0 & 0 & x\\
0 & 0 & x & 0 & 0 & 1 & 0 & 0\\
0 & 0 & 0 & 0 & 0 & 0 & 1 & 0\\
0 & 0 & 0 & x & 0 & 0 & 0 & 1\\
\end{pmatrix}
\,
M'_{7}=
\begin{pmatrix}
1 & x & 0 & 0 & 0 & 0 & 0\\
x & 1 & 0 & 0 & 0 & 0 & 0\\
0 & 0 & 1 & 0 & x & 0 & 0\\
0 & x & 0 & 1 & 0 & 0 & 0\\
0 & 0 & x & 0 & 0 & 1 & 0\\
0 & 0 & 0 & 0 & 0 & 0 & 1\\
0 & 0 & 0 & x & 0 & 0 & 0\\
\end{pmatrix}
\,
\widetilde{M}_{7}=
\begin{pmatrix}
1 & x & 0 & 0 & 0 & 0 & 0\\
x & 1 & 0 & 0 & 0 & 0 & 0\\
0 & 0 & 1 & 0 & x & 0 & 0\\
0 & x & 0 & 1 & 0 & 0 & 0\\
0 & 0 & 0 & x & 0 & 0 & 0\\
0 & 0 & x & 0 & 0 & 1 & 0\\
0 & 0 & 0 & 0 & 0 & 0 & 1
\end{pmatrix}
\end{equation*}
 Proving that $\det M'_{k-1} = 0$ is equivalent to proving that $\det\widetilde{M}_{k-1}=0$. Moreover, the $\widetilde{M}_{k-1}$ is basically a Collatz matrix but for row $(2k-1)/3$. In order to prove that $\det\widetilde{M}_{k-1}=0$ we argue by contradiction. Along row $(2k-1)/3$, $\widetilde{M}_{k-1}$ has only one non-zero entry (look at the rightmost matrix above for $k=8$, where $(2k-1)/3=5$ and $k/2=4$). This is an $x$ along column $k/2$.  We recall that, by definition, $$\det\widetilde{M}_{k-1} = \underset{\tau \text{ is perm of } \{ 1, 2, \dots , k-1 \}}{\sum} \quad \text{sign} (\tau)\prod_{i=1}^{k-1}\big(\widetilde{M}_{k-1}\big)_{i, \tau (i)} $$ So, if $\det\widetilde{M}_{k-1}\neq 0$, there is a permutation of $\{ 1, 2, \dots , k-1 \}$ (call it $\sigma$) such that the corresponding summand   $$\text{sign} (\sigma)\prod_{i=1}^{k-1}\big(\widetilde{M}_{k-1}\big)_{i, \sigma (i)} $$  in the formula for the determinant of $\widetilde{M}_{k-1}$, is non-zero. Thus $$\sigma\bigg(\frac{2k-1}{3}\bigg) = \frac{k}{2}.$$ We now try to find the cycle of $\sigma^{-1}$ which contains $k/2$. The next element is $(2k-1)/3$. We hope to reach an absurd statement like such a cycle cannot exist. In passing, we are dealing with $\sigma^{-1}$ for there seems to be less branching of the possibilities than when dealing with $\sigma$. This methodology will not clear the remaining two cases ($k=44+54l$ and $k=26+54l$), but will clear subsequences of these numbers and hopefully will provide inspiration for advances on the Collatz conjecture. Also, it gives a simpler answer than that of \cite{JFA}'s for $k=8+54l$. We use the following terminology: $$k_{-1}:=\sigma^{-1}(k/2)(=(2k-1)/3), \qquad \qquad k_{-j}:=\sigma^{-j}(k/2)\qquad \text{where $\sigma^{-j}$ is the $j$-th iterate of $\sigma^{-1}$}.$$  Since the rows of $\widetilde{M}_{k-1}$ comply with the definition of rows of a Collatz matrix (but for row $(2k-1)/3$) the iterates of a positive integer $x\neq k/2$ will be provided by the maps $2x$ or $(2x-1)/3$: $$\sigma^{-1}(x)= 2x \qquad \qquad \text{ or } \qquad \qquad \sigma^{-1}(x)=\frac{2x-1}{3} .$$ Sometimes both possibilities will be acceptable. We note that we must check that $$\sigma^{-1}(x)=2x < k \qquad \qquad  \text{ and that } \qquad \qquad \sigma^{-1}(x)=\frac{2x-1}{3} \quad \text{ yields an integer.} $$ Moreover, we must check if $\sigma^{-1}(x)=k/2$ for some $x$; if this occurs, it means there is a cycle and our argument by contradiction will not work.
\bigbreak

We use this methodology to prove (again) the case $k=8+54l$ before we use it to present new results. We believe our notation below is straightforward. Namely, arrows with symbols above them ending with a question mark indicate we are trying an inverse function on the obvious argument; an ${\bf X}$ at the end means this inverse does not work.  An ``inverse does not work'' when its image is larger than (or equal to) $k$ or its image is not an integer, according to the discussion above. If upon trying both inverses ($2x$ and $(2x-1)/3$) the ``inverse does not work'', we know this candidate to a cycle cannot reach its initial term ($k/2$). We thus conclude that such a cycle cannot exist and therefore the $\sigma$ cannot correspond to a non-null summand in the formula for the determinant. The final conclusion is that such determinant has to be zero.

\subsection{The case $k=8+54l$ (alternative method to that of \cite{JFA}).}\label{subsec:k=8+54l}
\bigbreak

\boxed{k = 8+54l, \qquad k_{-1}=5+36l, \qquad k/2=4+27l}

\bigbreak
\begin{align*}
& k_{-1}=\mathbf{5+36l} \qquad \xrightarrow{\text{$\times 2?$}}10+72l>8+54l \qquad {\bf X}\\
&\\
& k_{-2}=\frac{9+72l}{3}=\mathbf{3+24l}\qquad \xrightarrow{\text{$\times 2 - 1 / 3 ?$}} \frac{5+48l}{3}=\frac{5}{3}+16l\quad {\bf X}\\
&\\
& k_{-3}={\bf 6+48l} \qquad \xrightarrow{\text{$\times 2?$}}12+96l > 8+54l \qquad {\bf X}\\
&\\
& k_{-4}=\frac{11+96l}{3}=\frac{11}{3}+32l \qquad  {\bf X}
\end{align*}
\bigbreak
We thus conclude that these iterates do not form a cycle. Hence, $$\sigma\bigg(\frac{2k-1}{3}\bigg)\neq \frac{k}{2} .$$

Furthermore, we proved that, for $k=8+54l$, any permutation in the formula for the determinant of $\det\widetilde{M}_{k-1}$  satisfies $\sigma\big(\frac{2k-1}{3}\big)\neq \frac{k}{2}$. Since the only non-null entry in row $\frac{2k-1}{3}$ is along column $\frac{k}{2}$, this implies

$$\det \widetilde{M}_{k-1}=0 \qquad \qquad \qquad \qquad \text{ for } k=8+54l .$$  Hence,  $$\det M_{k}=\det M_{k-1} \qquad \qquad \text{ for } k=8+54l .$$ Again, this situation had already been cleared in \cite{JFA}. They solved it by proving that the linear homogeneous system of equations whose coefficient matrix is $\widetilde{M}_{k-1}$  (for $k=8+54l$) has a non-trivial solution. We think our method provides a simpler solution and applied to this case ($k=8+54l$), paves the way to showing how our argument works.

\bigbreak

By this analysis we see that to show that $\det\widetilde{M}_{k-1}=0$ it is sufficient to show that k/2 does not lie on a specific closed (inverse) Collatz orbit. In the calculations below we show directly (by orbit analysis) that this is the issue for the cases described in our Theorem \ref{thm:results}.
\subsection{The case $k=44+54l$.}\label{subsec:k=44+54l}
Now for the case $k=44+54l$. This case was not dealt with before.

\bigbreak
\boxed{k = 44+54l, \qquad k_{-1}=29+36l, \qquad k/2=22+27l}

\bigbreak
\begin{align*}
& k_{-1}=\mathbf{29+36l} \qquad \xrightarrow{\text{$\times 2?$}}58+72l>44+54l \qquad {\bf X}\\
&\\
& k_{-2}=\frac{57+72l}{3}=\mathbf{19+24l}\qquad \xrightarrow{\text{$\times 2 - 1 / 3 ?$}} \frac{37+48l}{3}=\frac{37}{3}+16l\quad {\bf X}\\
&\\
& k_{-3}={\bf 38+48l} \qquad \xrightarrow{\text{$\times 2?$}}76+96l > 44+54l \qquad {\bf X}\\
&\\
& k_{-4}=\frac{75+96l}{3}=\mathbf{25+32l} \qquad \xrightarrow{\text{$\times 2?$}}50+64l > 44+54l \quad {\bf X}\\
&\\
& k_{-5}=\frac{49+64l}{3}=\frac{177+64(l-2)}{3}=59+64\frac{l-2}{3}  \qquad {\bf X} \text{ unless } \boxed{l=2+3l_1}\\
\end{align*}

\bigbreak
\bigbreak
So, at this point we can state:
\bigbreak
\bigbreak
\boxed{\mathbf{ \text{For }k=44+54l, \text{ if } 3\, |\, l \text{ or } 3\, |\, (l-1),  \text{ then }  \det M_k = \det M_{k-1}.}}
\bigbreak
This is statement $(a)(i)$ in Theorem \ref{thm:results}.
\bigbreak
We now explore further the other case:
\bigbreak
\bigbreak
Update:
\bigbreak
\boxed{l=2+3l_1}
\bigbreak

so

\begin{align*}
&k = 44+54l=44+54(2+3l_1)=152+162l_1,\\
&k_{-1}=29+36l=29+36(2+3l_1)=101+108l_1,\\
&k/2=76+81l_1\\
\end{align*}

\begin{align*}
&k_{-5} = 59+64l_1\\
&\\
& k_{-6}^1=2[59+64l_1]=118+128l_1 \qquad \text{ or } \qquad k_{-6}^2=\frac{117+128l_1}{3}=39+128\frac{l_1}{3}\quad{\bf X}\quad \text{ unless } \boxed{l_1=3l_2^1}\\
&\\
& k_{-7}^1=\frac{235+256l_1}{3}=\frac{747+256(l_1-2)}{3}=249+256\frac{l_1-2}{3}\quad{\bf X}\quad \text{ unless } \boxed{l_1=2+3l_2^2} \\
\end{align*}

\bigbreak
\bigbreak
So, at this point we can state:
\bigbreak
\bigbreak
\boxed{\mathbf{ \text{For }k=44+54l, \text{ with } l=2+3l_1,  \text{ if } 3\,|\,(l_1-1),  \text{ then }  \det M_k = \det M_{k-1}.}}
\bigbreak
This is statement $(a)(ii)$ in Theorem \ref{thm:results}.
\bigbreak
We now explore further the other cases: $$l_1=3l_2^1\qquad \text{ or } \qquad l_1=2+3l_2^2 .$$
\bigbreak
\bigbreak

Update $1$: \qquad \qquad \qquad \qquad \qquad \qquad \qquad \qquad  \qquad Update $2$:
\bigbreak
\boxed{l_1=3l_2^1} \qquad \qquad \qquad \qquad  \qquad  \qquad \qquad \qquad \qquad  \quad \boxed{l_1=2+3l_2^2}
\bigbreak

so

\begin{align*}
&k = 152+162l_1 = 152+162(3l_2^1)=152+486l_2^1, \quad \quad  k = 152+162l_1 = 152+162(2+3l_2^2)=476+486l_2^2,\\
&k_{-1}=101+108(3l_2^1)=101+324l_2^1, \qquad \qquad \qquad \quad k_{-1}=101+108(2+3l_2^2)=317+324l_2^2,\\
&k/2=76+243l_2^1, \qquad \qquad \qquad \qquad \qquad \qquad \qquad \quad k/2=238+243l_2^2\\
\end{align*}

\bigbreak
\bigbreak

For Update $1$:

\bigbreak
\bigbreak

\begin{align*}
&k_{-6}^2 = 39+128l_2^1\\
&\\
& k_{-7}^{21}=78+256l_2^1 \qquad \text{ or } \qquad k_{-7}^{22}=\frac{77+256l_2^1}{3}=111+256\frac{l_2^1-1}{3}\quad{\bf X}\quad \text{ unless } \boxed{l_2^1=1+3l_3^2}\\
&\\
& k_{-8}^{21}=\frac{155+512l_2^1}{3}=393+512\frac{l_2^1-2}{3}\quad{\bf X}\quad \text{ unless } \boxed{l_2^1=2+3l_3^1} \\
\end{align*}

\bigbreak
\bigbreak
So, at this point we can state:
\bigbreak
\bigbreak
\boxed{\mathbf{ \text{For }k=44+54l, \text{ with } l=2+3l_1,  l_1=3l_2^1\, \text{ if } 3\,|\,l_2^1,  \text{ then }  \det M_k = \det M_{k-1}.}}
\bigbreak
This is statement $(a)(iii)(1)$ in Theorem \ref{thm:results}.
\bigbreak

For Update $2$:

\bigbreak
\bigbreak

\begin{align*}
&k_{-7}^1 = 249+256l_2^2\\
&\\
& k_{-8}^{1}=\frac{497+512l_2^2}{3}=507+512\frac{l_2^2-2}{3}\quad{\bf X}\quad \text{ unless } \boxed{l_2^2=2+3l_3^3}\\
\end{align*}
\bigbreak
\bigbreak
So, at this point we can state:
\bigbreak
\bigbreak
\boxed{\mathbf{ \text{For }k=44+54l, \text{ with } l=2+3l_1,  l_1=2+3l_2^2\, \text{ if } 3\,|\,l_2^2,  \text{ or } 3\,|\,(l_2^2-1),  \text{ then }  \det M_k = \det M_{k-1}.}}
\bigbreak
This is statement $(a)(iii)(2)$ in Theorem \ref{thm:results}.
\bigbreak

\subsection{The case $k=26+54l$.}\label{subsec:k=26+54l}

Now for the case $k=26+54l$. This case was not dealt with before.

\bigbreak
\boxed{k = 26+54l, \qquad k_{-1}=17+36l, \qquad k/2=13+27l}

\bigbreak
\begin{align*}
& k_{-1}=\mathbf{17+36l} \qquad \xrightarrow{\text{$\times 2?$}}34+72l>26+54l \qquad {\bf X}\\
&\\
& k_{-2}=\frac{33+72l}{3}=\mathbf{11+24l}\\
&\\
& k_{-3}^1=\mathbf{22+48l} \qquad \text{ or } \qquad k_{-3}^2=\frac{21+48l}{3}=\mathbf{7+16l}\\
&\\
& k_{-4}^1=\frac{43}{3}+32l\qquad {\bf X}\\
&\\
& k_{-4}^{21}=14+32l \qquad \text{ or } \qquad k_{-4}^{22}=\frac{13+32l}{3}=15+32\frac{l-1}{3}  \qquad {\bf X} \text{ unless } \boxed{l=1+3l_1^1}\\
&\\
& k_{-5}^{21}=\frac{27+64l}{3}=9+64\frac{l}{3} \qquad {\bf X} \text{ unless } \boxed{l=3l_1^2}\\
\end{align*}

\bigbreak
\bigbreak
So, at this point we can state:
\bigbreak
\bigbreak
\boxed{\mathbf{ \text{For }k=26+54l, \text{ if }  3\, |\, (l-2),  \text{ then }  \det M_k = \det M_{k-1}.}}
\bigbreak
This is statement $(b)(i)$ in Theorem \ref{thm:results}.
\bigbreak
\bigbreak
\bigbreak
We now explore further the other cases: $$l=1+3l_1^1\qquad \text{ or } \qquad l=3l_1^2 .$$
\bigbreak
\bigbreak

Update $1$: \qquad \qquad \qquad \qquad \qquad \qquad \qquad \qquad  \qquad Update $2$:
\bigbreak
\boxed{l=1+3l_1^1} \qquad \qquad \qquad \qquad  \qquad  \qquad \qquad \qquad \qquad  \quad \boxed{l=3l_1^2}
\bigbreak

so

\begin{align*}
&k = 26+54l = 26+54(1+3l_1^1)=80+162l_1^1, \qquad \quad  k = 26+54l = 26+54(3l_1^2)=26+162l_1^2,\\
&k_{-1}=53+108l_1^1, \qquad \qquad \qquad \qquad \qquad \qquad \qquad \quad k_{-1}=17+108l_1^2,\\
&k/2=40+81l_1^1, \qquad \qquad \qquad \qquad \qquad \qquad \qquad \qquad k/2=13+81l_1^2\\
\end{align*}

\bigbreak
\bigbreak

For Update $1$:

\bigbreak
\bigbreak
\begin{align*}
&k_{-4}^{22} = 15+32l_1^1\\
&\\
& k_{-5}^{221}=30+64l_1^1 \qquad \text{ or } \qquad k_{-5}^{222}=\frac{29+64l_1^1}{3}=31+64\frac{l_1^1-1}{3}\quad{\bf X}\quad \text{ unless } \boxed{l_1^1=1+3l_3^1}\\
&\\
& k_{-6}^{2211}=60+128l_1^1 \qquad \text{ or } \qquad k_{-6}^{2212}=\frac{59+128l_1^1}{3}=105+128\frac{l_1^1-2}{3}\quad{\bf X}\quad \text{ unless } \boxed{l_1^1=2+3l_3^2} \\
& k_{-7}^{2211}=\frac{119+256l_1^1}{3}=125+256\frac{l_1^1-1}{3}\quad{\bf X}\quad \text{ unless } \boxed{l_1^1=1+3l_3^3} \\
\end{align*}

\bigbreak
\bigbreak
So, at this point we can state:
\bigbreak
\bigbreak
\boxed{\mathbf{ \text{For }k=26+54l, \text{ with } l=1+3l_1^1,   \text{ if } 3\,|\,l_1^1,  \text{ then }  \det M_k = \det M_{k-1}.}}
\bigbreak
This is statement $(b)(ii)(1)$ in Theorem \ref{thm:results}.
\bigbreak

For Update $2$:

\bigbreak
\bigbreak
\begin{align*}
&k_{-5} = 9+64l_1^2\\
&\\
& k_{-6}^{1}=18+128l_1^2 \qquad \text{ or } \qquad k_{-6}^{2}=\frac{17+128l_1^2}{3}=91+128\frac{l_1^2-2}{3}\quad{\bf X}\quad \text{ unless } \boxed{l_1^2=2+3l_3^2}\\
&\\
& k_{-7}^{11}=\frac{35+256l_1^2}{3}=97+256\frac{l_1^2-1}{3}\quad{\bf X}\quad \text{ unless } \boxed{l_1^2=1+3l_3^1} \\
\end{align*}

\bigbreak
\bigbreak
So, at this point we can state:
\bigbreak
\bigbreak
\boxed{\mathbf{ \text{For }k=26+54l, \text{ with } l=3l_1^2,   \text{ if } 3\,|\,l_1^2,  \text{ then }  \det M_k = \det M_{k-1}.}}
\bigbreak
This is statement $(b)(ii)(2)$ in Theorem \ref{thm:results}.
\bigbreak

\section{Final Remarks}\label{sec:Final}

We hope that by a deeper look into this matrix reformulation, the Collatz Orbit Conjecture will be fully resolved.

\end{document}